\documentclass[12pt]{article}

\usepackage[dvips]{graphicx}
\usepackage{amsmath,amssymb}
\usepackage{subfigure}
\usepackage{natbib}

\newtheorem{theorem}{Theorem}[section]
\newtheorem{lemma}[theorem]{Lemma}
\newtheorem{corollary}[theorem]{Corollary}
\newtheorem{proposition}[theorem]{Proposition}

\newtheorem{example}[theorem]{Example}

\newcommand{\qed}{ \hfill $\square$ \\ }

\newcommand{\diff}{\mathrm{d}}

\def\Proof{\noindent{\bf Proof.}\quad}

\begin{document}
\title{Bayesian shrinkage prediction for the regression problem}

\author{Kei Kobayashi\footnote{kei@ism.ac.jp}
and Fumiyasu Komaki\footnote{komaki@mist.i.u-tokyo.ac.jp}
}

\date{}

\maketitle

\begin{abstract}
We consider Bayesian shrinkage predictions
for the Normal regression problem under the frequentist Kullback-Leibler risk function.

Firstly,  we consider the multivariate Normal model with an unknown mean and
a known covariance.
While the unknown mean is fixed, the covariance of future samples can be
different from training samples.
We show that the Bayesian predictive distribution
based on the uniform prior is dominated by that based on a class of priors
if the prior distributions for the covariance and future covariance matrices
are rotation invariant.

Then, we consider a class of priors for the mean parameters depending on
the future covariance matrix. 
With such a prior, we can construct a Bayesian predictive distribution dominating
that based on the uniform prior.

Lastly, applying this result to the prediction
of response variables in the Normal linear regression model, 
we show that there exists a Bayesian predictive distribution
dominating that based on the uniform prior.
Minimaxity of these Bayesian predictions follows from these results.

\end{abstract}

\vspace{3mm}

\section{Introduction}
Suppose that we have observations $y\sim N_d(y;\mu,\Sigma)$.
Here $N_d$ is the density function of the
$d$-dimensional multivariate Normal distribution with mean vector $\mu$ and 
covariance matrix $\Sigma$.
We consider the prediction of $\tilde{y}\sim N_d(\tilde{y};\mu,\tilde{\Sigma})$
using a predictive density $\hat{p}(\tilde{y} | y).$
We assume that the mean of the distribution of unobserved (future) samples is
the same as the one of the observed samples.
However, the covariance matrices, $\Sigma$ and $\tilde{\Sigma}$, are
not necessarily the same or proportional to each other.
We call a problem with such settings the ``problem with changeable covariances.''
As we will show below, the changeable covariance is a natural assumption
when we consider the linear regression problems.

In the present work, we assume that the mean vector $\mu$ is unknown
and the covariance matrix $\Sigma$ is known.
We consider both cases where the future covariance $\tilde{\Sigma}$ is known
and unknown.

We evaluate predictive densities $\hat{p}(\tilde{y} | y)$
by the KL loss function
\begin{equation}
D(\tilde{p}(\tilde{y}|\theta) \| \hat{p}(\tilde{y}|y))
:=
\int \tilde{p}(\tilde{y}|\theta) \log \frac{\tilde{p}(\tilde{y}|\theta)}
{\hat{p}(\tilde{y}|y)} \diff \tilde{y}
\label{eq:reg-KL-loss}
\end{equation}
and the (frequentist) risk function 
\begin{equation}
R_{\rm KL}(\hat{p},\theta)
:= \int p(y|\theta) D(\tilde{p}(\tilde{y}|\theta) \|
 \hat{p}(\tilde{y}|y)) \diff \tilde{y}.
\label{eq:reg-KL-risk}
\end{equation}

We consider the Bayesian predictive density
\begin{equation*}
p_\pi(\tilde{y}|y):=
\frac{\int 
\tilde{p}(\tilde{y}|\theta)p(y|\theta)\pi(\theta)\diff \theta}
{\int p(y|\theta)\pi(\theta)\diff \theta}
\end{equation*}
with prior $\pi(\theta)$.
For the Normal model,
the Bayesian predictive density with the uniform prior $\pi_I(\mu)=1$ becomes
\begin{equation*}
p_\pi(\tilde{y}|y;\Sigma,\tilde{\Sigma})=
 \frac{1}{(2\pi)^{d/2}|\Sigma+\tilde{\Sigma}|^{1/2}}
\exp\Big(-\frac{(\tilde{y}-y)^\top (\Sigma+\tilde{\Sigma})^{-1}(\tilde{y}-y)}
{2}  \Big),
\end{equation*}
as we will see in Section \ref{sec:main_result}.
Let $p_\pi(\tilde{y}|y;\Sigma,\tilde{\Sigma})$ denote $p_\pi(\tilde{y}|y)$
for short.

When $\tilde{\Sigma}$ is proportional to $\Sigma$, i.e. $\tilde{\Sigma}=a\Sigma$
for $a>0$, the problem is reduced to the one with $\Sigma=v I_d$ and 
$\tilde{\Sigma}= \tilde{v} I_d$ for positive scalar values $v$ and $\tilde{v}$. 
This case with `unchangeable covariances' has been well studied.
The Bayesian predictive density
\begin{equation*}
p_I(\tilde{y}|y;\Sigma,\tilde{\Sigma})=
 \frac{1}{\{2\pi(v+\tilde{v})\}^{d/2}}
\exp\Big(-\frac{\|\tilde{y}-y\|^2}
{2(v+\tilde{v})}  \Big)
\end{equation*}
based on the uniform prior $\pi_{\rm I}(\mu)=1$
dominates the plug-in density
$$
p(\tilde{y}|\hat{\mu})= \frac{1}{\{2\pi \tilde{v}\}^{d/2}}
                        \exp\Big(-\frac{\|\tilde{y}-y\|^2}{2\tilde{v}} \Big)
$$
with MLE,
where $\hat{\mu}=y$.
Moreover, by \cite{murray1977} and \cite{ng1980},
the Bayesian predictive density
$p_{\rm I}(\tilde{y}|y)$ is the best predictive density that
is invariant under the translation group.
In \cite{liang-barron-2004IEEEIT} and \cite{george_etc2006},
the minimaxity of $p_{\rm I}$ was proved.

In \cite{komaki2001}, it was proved that the Bayesian predictive density
$p_{\rm S}(\tilde{y}|y)$ with Stein prior
\begin{equation}
\pi_{\rm S}(\mu):= \|\mu\|^{-(d-2)}
\label{eq:reg-Stein-prior}
\end{equation}
dominates the Bayesian predictive density $p_{\rm I}(\tilde{y}|y)$ with the
uniform prior $\pi_{\rm I}(\mu)$.

\cite{george_etc2006} generalized the result of \cite{komaki2001}.
Define the marginal distribution $m_\pi$ by
\begin{equation}
m_\pi(z;\Sigma) := \int {\rm N}(z;\mu, \Sigma)\pi(\mu) \,\diff\mu .
\label{eq:reg-marginal-dist}
\end{equation}
As we will see in Theorem \ref{thm:reg-1} below,
\cite{george_etc2006} proved a sufficient 
condition on the prior $\pi(\mu)$
or the marginal distribution $m_\pi$
for $p_\pi(\tilde{y}|y)$ to dominate $p_{\rm I}(\tilde{y}|y)$
when $\Sigma$ is proportional to $\tilde{\Sigma}$. 
In the present work,  we generalize the results
of \cite{komaki2001} and \cite{george_etc2006}
to the corresponding problem with the changeable covariances,
considering only finite sample cases.
Asymptotic properties of Bayesian prediction are studied
in \cite{komaki1996}, \cite{corcuera_giummole2000},
and \cite{komaki-2006AS}.

\section{Prior distributions independent of the future covariance}
\label{sec:main_result}
In this section, we develop and prove our main results
concerning properties of $p_\pi(\tilde{y}|y)$ in the
problem with changeable covariances.

First we give three lemmas generalizing
results proved in \cite{george_etc2006}
for the problem with ``unchangeable'' variances.

Define the marginal distribution $m_\pi$ by (\ref{eq:reg-marginal-dist}).

\begin{lemma}
\label{lem:reg-1}
If $m_\pi(z;\Sigma)<\infty$ for all $z$, then 
$p_\pi(\tilde{y}|y)$ is a proper probability density.
Moreover, the mean of $p_\pi(\tilde{y}|y)$ is equal to the posterior mean
$E_\pi[\mu | y]$ if it exists.
\end{lemma}

Let
$$w:=(\Sigma^{-1}+\tilde{\Sigma}^{-1})^{-1}
(\Sigma^{-1}y+\tilde{\Sigma}^{-1}\tilde{y})$$
and
\begin{equation}
\Sigma_w := (\Sigma^{-1}+\tilde{\Sigma}^{-1})^{-1}.
\label{eq:Sig_w}
\end{equation}

As a function of the predictive density
based on the uniform prior,
the Bayesian predictive density
based on a prior $\pi(\mu)$ becomes as follows:
\begin{lemma}
\label{lem:reg-2}
$$
p_\pi(\tilde{y}|y)= p_{\rm I}(\tilde{y}|y)
\frac{m_{\pi}(w;\Sigma_w)}{m_\pi(y;\Sigma)}.
$$
\end{lemma}

The following lemma is used for proving minimaxity of $p_{\pi}(\tilde{y}|y)$.
\begin{lemma}
\label{lem:reg-3}
The Bayesian predictive density $p_I(\tilde{y}|y)$ is minimax
under KL risk function $R_{\rm KL}(\hat{p},\mu)$.
\end{lemma}

Since the proofs of Lemma \ref{lem:reg-1} and Lemma \ref{lem:reg-3}
are almost same as those of Lemma 1 and Lemma 3 in \cite{george_etc2006},
we omit them. We prove only Lemma \ref{lem:reg-2}.

\quad\\
\noindent{\bf Proof of Lemma \ref{lem:reg-2}}\quad
\begin{align}
p(y|&\mu,\Sigma) p(\tilde{y}|\mu,\tilde{\Sigma})\nonumber\\
&=\frac{1}{(2\pi)^{d/2}|\Sigma|^{1/2}}
\exp\Big(-\frac{(y-\mu)^{\top} \Sigma^{-1} (y-\mu)}{2}\Big)
\frac{1}{(2\pi)^{d/2}|\tilde{\Sigma}|^{1/2}}
\exp\Big(-\frac{(\tilde{y}-\mu)^{\top} \tilde{\Sigma}^{-1} (\tilde{y}-\mu)}{2}\Big)\nonumber\\
&=\frac{1}{(2\pi)^{d/2}|\Sigma|^{1/2}}\frac{1}{(2\pi)^{d/2}|\tilde{\Sigma}|^{1/2}}
\exp\Big(-\frac{(w-\mu)^{\top}\Sigma_w^{-1}(w-\mu)}{2}\Big)
\exp\Big(-\frac{y^{\top}\Sigma^{-1}y}{2}-\frac{\tilde{y}^{\top}\tilde{\Sigma}^{-1}\tilde{y}}{2}\Big)
\nonumber\\
&~~~~~~~~~~~~~~~~~~~~~~~~~~~~~~~~~~~~~~~~~~~~~~~
\exp\Big(\frac{(\Sigma^{-1}y+\tilde{\Sigma}^{-1}\tilde{y})^{\top}
(\Sigma^{-1}+\tilde{\Sigma}^{-1})^{-1}(\Sigma^{-1}y+\tilde{\Sigma}^{-1}\tilde{y})
}{2}\Big)\nonumber\\
&=\frac{1}{(2\pi)^{d/2}|\Sigma|^{1/2}}\frac{1}{(2\pi)^{d/2}|\tilde{\Sigma}|^{1/2}}
\exp\Big(-\frac{(w-\mu)^{\top}\Sigma_w^{-1}(w-\mu)}{2}\Big)
\exp\Big(-\frac{(y-\tilde{y})^{\top}(\Sigma+\tilde{\Sigma})^{-1}(y-\tilde{y})}{2}\Big)
\label{eq_lem1-1}
\end{align}
In the last equation, we use
\begin{align*}
\Sigma^{-1}&(\Sigma^{-1}+\tilde{\Sigma}^{-1})^{-1}\Sigma^{-1}-\Sigma^{-1}\\
&=\Sigma^{-1}(\Sigma^{-1}+\tilde{\Sigma}^{-1})^{-1}\Sigma^{-1}
-\Sigma^{-1}(\Sigma^{-1}+\tilde{\Sigma}^{-1})^{-1}(\Sigma^{-1}+\tilde{\Sigma}^{-1})\\
&=-\Sigma^{-1}(\Sigma^{-1}+\tilde{\Sigma}^{-1})^{-1}\tilde{\Sigma}^{-1}\\
&=-(\Sigma+\tilde{\Sigma})^{-1}.
\end{align*}

From (\ref{eq_lem1-1}), the predictive density with the uniform prior 
$I(\mu)=1$ is given by 
\begin{align*}
p_{\rm I}(\tilde{y}|y)&= \frac{\int p(y|\mu,\Sigma) p(\tilde{y}|\mu,\tilde{\Sigma})\,\diff\mu}
{\int p(y|\mu, \Sigma)\,\diff\mu}\\
&=\frac{1}{(2\pi)^{d/2}|\Sigma|^{1/2}}\frac{1}{(2\pi)^{d/2}|\tilde{\Sigma}|^{1/2}}
|\Sigma^{-1}+\tilde{\Sigma}^{-1}|^{-1/2}(2\pi)^{d/2} 
\exp\Big(-\frac{(y-\tilde{y})^{\top}(\Sigma+\tilde{\Sigma})^{-1}(y-\tilde{y})}{2}\Big)\\
&=(2\pi)^{-d/2}|\Sigma+\tilde{\Sigma}|^{-1/2}
\exp\Big(-\frac{(y-\tilde{y})^{\top}(\Sigma+\tilde{\Sigma})^{-1}(y-\tilde{y})}{2}\Big).
\end{align*}

Therefore
\begin{align*}
p_{\pi}(\tilde{y}|y)&= \frac{\int p(y|\mu,\Sigma) p(\tilde{y}|\mu,\tilde{\Sigma})\pi(\mu)
\,\diff\mu}
{\int p(y|\mu, \Sigma)\pi(\mu) \,\diff\mu}\\
&=\frac{p_{\rm I}(\tilde{y}|y)\int {\rm N}(w;\mu,\Sigma_w)\pi(\mu)\,\diff\mu}
{\int {\rm N}(y;\mu,\Sigma)\pi(\mu) \,\diff\mu}\\
&=p_{\rm I}(\tilde{y}|y)\frac{m_{\pi}(w;\Sigma_w)}{m_\pi(y;\Sigma)}.
\end{align*}
\qed

Next, the difference of the risk functions of the two priors is evaluated.
Let
$$
R_{\rm KL}(\pi,\mu):= \int p(y|\mu,\Sigma)D(p(\tilde{y}|\mu,\tilde{\Sigma})
\|  p_\pi(\tilde{y}|y)) \,\diff y
$$
$$
\phi_\pi(\mu,\Sigma) := \int {\rm N}(z;\mu,\Sigma)\log m_\pi(z;\Sigma) 
\,\diff z.
$$
Then from Lemma \ref{lem:reg-2},
\begin{align}
R_{\rm KL}(\pi,\mu)-R_{\rm KL}(\pi_{\rm I},\mu)
&=\int p(y|\mu,\Sigma)p(\tilde{y}|\mu,\tilde{\Sigma})\log 
\frac{p_{\rm I}(\tilde{y}|y)}{p_\pi(\tilde{y}|y)} \,\diff y \,\diff \tilde{y}
\nonumber\\
&=\int p(y|\mu,\Sigma)p(\tilde{y}|\mu,\tilde{\Sigma})\log 
\frac{m_{\pi}(y;\Sigma)}{m_\pi(w;\Sigma_w)} \,\diff y \,\diff \tilde{y}
\nonumber\\
&=\phi_\pi(\mu,\Sigma)-\phi_\pi(\mu,\Sigma_w).
\label{eq:reg-risk-difference}
\end{align}
Now $\Sigma_{w}=(\Sigma^{-1}+\tilde{\Sigma}^{-1})^{-1} \prec \Sigma$.
In order to prove $R_{\rm KL}(\pi,\mu)< R_{\rm KL}(\pi_{\rm I},\mu)$,
it suffices to prove $\phi_\pi(\mu,\Sigma)< \phi_\pi(\mu,\Sigma_w)$.

Before stating the main results for the problem with changeable covariances,
we review some results with a special setting, i.e., unchangeable covariances.

An extended real-valued function $\pi(\mu)$ on an open set 
$R\subset \mathbb{R}^p$ is said to be {\it superharmonic}
when it satisfies the following properties:
\begin{enumerate}
\item 
$-\infty<\pi(\mu)\leq \infty$ and $\pi(\mu)\not\equiv \infty$ on
any component of $R$.
\item
$\pi(\mu)$ is lower semi-continuous on $R$.
\item
If $G$ is an open subset of $R$ with compact closure $\Bar{G}\subset R$,
$w(\mu)$ is a continuous function on $\Bar{G}$, $w(\mu)$ is harmonic on $G$,
and $\pi(\mu)\geq w(\mu)$ on $\partial{G}$, then $\pi(\mu)\geq w(\mu)$
on $G$.
\end{enumerate}

If $\pi(\mu)$ is a $C^2$ function, then $\pi(\mu)$
is superharmonic on $R$ if and only if $\Delta \pi\leq 0$ on $R$.

\begin{theorem}[\cite{komaki2001} and \cite{george_etc2006}]
\label{thm:reg-1}\quad

Assume $d\geq 3$.

(i) 
If $\pi(\mu)$ is the Stein prior $\pi_{\rm S}(\mu)$,
$$
v_1 > v_2 > 0
\Rightarrow
\phi_\pi(\mu,v_1 I_d) < \phi_\pi(\mu,v_2 I_d) \mbox{ for all } \mu.
$$

(ii) 
If $\pi(\mu)$ is a superharmonic function
and $m_\pi(z;v I_d)<\infty$ for any $z$ and $v$, 
$$
v_1 > v_2 > 0
\Rightarrow
\phi_\pi(\mu,v_1 I_d) \leq \phi_\pi(\mu,v_2 I_d) \mbox{ for all } \mu.
$$
~~~~~~~~Furthermore, if $m_\pi(z;v I_d)$ is also not constant for all
$v_2 \leq v \leq v_1$, the inequality\\
~~~~~~~~holds strictly.

(iii)
If $\sqrt{m_\pi(z;v I_d)}$ is a superharmonic function for any $v$ and 
$m_\pi(z;v I_d)<\infty$ for any $z$ and $v$, 
$$
v_1 > v_2 > 0
\Rightarrow
\phi_\pi(\mu,v_1 I_d) \leq \phi_\pi(\mu,v_2 I_d) \mbox{ for all } \mu.
$$
~~~~~~~~Furthermore, if $m_\pi(z;v I_d)$ is also not constant for any
$v_2 \leq v \leq v_1$, the inequality\\
~~~~~~~~holds strictly.
\end{theorem}

We note that (iii) implies (ii) and (ii) implies (i).
(i) was proved in \cite{komaki2001}.
(ii) and (iii) were proved in \cite{george_etc2006}.

Theorem \ref{thm:reg-rot-2} is a generalization of
(ii) of Theorem \ref{thm:reg-1} to the problem
with changeable covariances.
For each prior $\pi(\mu)$, define {\it a rescaled prior} with 
respect to a positive definite $d\times d$ matrix $\Sigma^*$ by
$$\pi_{\Sigma^*}(\mu) := \pi(\Sigma^{* -1/2}\mu).$$
In particular, call
$\pi_{{\rm S}; \Sigma^*}(\mu):= \pi_{\rm S}(\Sigma^{* -1/2}\mu)$
as {\it a rescaled Stein prior} with respect to $\Sigma^*$.

We consider Bayesian risk with priors $p(\Sigma)$ and $\tilde{p}(\tilde{\Sigma})$:
$$
\mathcal{R}_{\rm KL}(\pi,\mu)=\int 
p(\Sigma)\tilde{p}(\tilde{\Sigma})R_{\rm KL}(\pi,\mu) \diff \Sigma \diff \tilde{\Sigma},
$$
where $\diff \Sigma$ means a Lebesgue measure for a vector space of all 
components of a matrix $\Sigma$.
Define
\begin{align}
\varphi_\pi(\mu) &:= \int p(\Sigma)\tilde{p}(\tilde{\Sigma})
\phi_\pi(\mu,\Sigma) \diff \Sigma \diff \tilde{\Sigma} \nonumber\\
&=
\int p(\Sigma)\tilde{p}(\tilde{\Sigma}) {\rm N}(z;\mu,\Sigma)\log m_\pi(z;\Sigma) 
\,\diff z \diff \Sigma \diff \tilde{\Sigma}
\end{align}
\begin{align}
\varphi_\pi^w(\mu) &:= \int p(\Sigma)\tilde{p}(\tilde{\Sigma})
\phi_\pi(\mu,\Sigma_w) \diff \Sigma \diff \tilde{\Sigma} \nonumber\\
&=
\int p(\Sigma)\tilde{p}(\tilde{\Sigma}) {\rm N}(z;\mu,\Sigma_w)\log m_\pi(z;\Sigma_w) 
\,\diff z \diff \Sigma \diff \tilde{\Sigma}.
\end{align}

Then from (\ref{eq:reg-risk-difference}),
\begin{align}
\mathcal{R}_{\rm KL}(\pi,\mu)-
\mathcal{R}_{\rm KL}(\pi_{\rm I},\mu)
=\varphi_\pi(\mu)-\varphi_\pi^w(\mu).
\label{eq:reg-bayes-risk-difference}
\end{align}

We consider the case where $p(\Sigma)$, $\tilde{p}(\tilde{\Sigma})$, and $\pi(\mu)$
are rotation invariant.
Here, a function $f(\Sigma)$ of a matrix $\Sigma\in \mathbb{R}^{d\times 
d}$ and a function $f(\mu)$ of a vector $\mu\in \mathbb{R}^{d\times d}$ 
are said to be {\it rotation invariant} if $f(\Sigma)=f(P\Sigma P^\top)$ 
and $g(\mu)=g(P\mu)$, respectively, for every orthogonal matrix $P\in 
\mathbb{R}^{d\times d}$.

\begin{theorem}~\\
\label{thm:reg-rot-2}
Let $d\geq 3$. If $p(\Sigma)$ and $\tilde{p}(\tilde{\Sigma})$
are rotation invariant functions and
$\pi$ is a rotation invariant superharmonic prior,
then 
$$\mathcal{R}_{\rm KL}(\pi_\Sigma,\mu)\leq \mathcal{R}_{\rm KL}(\pi_{\rm I},\mu)$$
for any $\mu$.
In particular, the Bayesian predictive distribution $p_\Sigma(y|\tilde{y})$
with $\pi_\Sigma$ dominates that based on $\pi_{\rm I}$ if $\pi$ is also not constant.

\end{theorem}

\Proof
We note that $m_\pi(z;\Sigma)<\infty$ for every $z\in \mathbb{R}^d$
and positive definite matrix $\Sigma \in \mathbb{R}^{d\times d}$ 
from Lemma \ref{lem:finite-marginal} in the appendix.

First, we prove invariance of 
$\varphi_{\pi_{\Sigma}}(\mu)$
and $\varphi_{\pi_{\Sigma}}^w(\mu)$
under rotations of $\mu$.

Let $P$ be a $d\times d$ orthogonal matrix, then
\begin{align*}
\varphi_{\pi_{\Sigma}}(P\mu) &= \int p(\Sigma)\tilde{p}(\tilde{\Sigma})
N(z;P\mu,\Sigma) \log \int N(z;\mu',\Sigma) \pi_{\Sigma}(\mu')
d\mu' dz d\Sigma d\tilde{\Sigma}
\nonumber\\
&=\int p(\Sigma)\tilde{p}(\tilde{\Sigma})
N(\tilde{z};\mu,P^\top \Sigma P) 
\log \int N(\tilde{z};\tilde{\mu}',P^\top \Sigma P)
\pi(\Sigma^{-1/2} P \tilde{\mu}')
d\tilde{\mu}' d\tilde{z} d\Sigma d\tilde{\Sigma}\\
&=
\int p(P \Sigma P^\top)\tilde{p}(\tilde{\Sigma})
N(\tilde{z};\mu,\Sigma) 
\log \int N(\tilde{z};\tilde{\mu}',\Sigma) \pi(\Sigma^{-1/2} \tilde{\mu}')
d\tilde{\mu}' d\tilde{z} d\Sigma d\tilde{\Sigma}\\
&= \varphi_{\pi_{\Sigma}}(\mu).
\end{align*}
Proof of the rotation invariance of $\varphi_{\pi_{\Sigma}}^w(\mu)$
is nearly the same.

We define
$$
\mu^* := \arg\max_{\|\mu'\|=\|\mu \|}
\frac{\|\Sigma^{-1/2}\mu'\|}{\|\Sigma_w^{-1/2}\mu'\|}
$$
and
\begin{equation}
\tau := \frac{\|\Sigma^{-1/2}\mu^*\|}{\|\Sigma_w^{-1/2}\mu^*\|}.
\label{eq:tau}
\end{equation}
Note that $0<\tau< 1$, because $\tilde{\Sigma}$ is positive definite.
Moreover,
$$\|\tau \Sigma_w^{-1/2}\tilde{\mu}'\|
= \tau \|\Sigma_w^{-1/2}\tilde{\mu}'\|
\geq \frac{\|\Sigma^{-1/2}\tilde{\mu}'\|}{\|\Sigma_w^{-1/2}\tilde{\mu}'\|}
\|\Sigma_w^{-1/2}\tilde{\mu}'\|= \|\Sigma^{-1/2}\tilde{\mu}'\|
$$ 
for every $\tilde{\mu}'$.

From the rotation invariance of $\phi_{\pi_\Sigma}$,
\begin{align}
\varphi_{\pi_{\Sigma}}(\mu) &= \varphi_{\pi_{\Sigma}}(\mu^*)\nonumber\\
&=E_{\Sigma,\tilde{\Sigma}}[\int N(z;\mu^*,\Sigma) \log \int N(z;\tilde{\mu},\Sigma)
\pi(\Sigma^{-1/2}\tilde{\mu})d\tilde{\mu} dz]\nonumber\\
&=E_{\Sigma,\tilde{\Sigma}}[\int N(\tilde{z};\Sigma^{-1/2}\mu^*,I_d) 
\log \int N(\tilde{z};\tilde{\mu}',I_d)\pi(\tilde{\mu}')
d\tilde{\mu}' d\tilde{z}]\nonumber\\
&=E_{\Sigma,\tilde{\Sigma}}[\int N(\tilde{z};\tau\Sigma_w^{-1/2}\mu^*,I_d) 
\log \int N(\tilde{z};\tilde{\mu}',I_d)\pi(\tilde{\mu}')
d\tilde{\mu}' d\tilde{z}]\nonumber\\
&=E_{\Sigma,\tilde{\Sigma}}[\int N(\tilde{z};\Sigma_w^{-1/2}\mu^*,\tau^{-2}I_d) 
\log \int N(\tilde{z};\tilde{\mu}',\tau^{-2}I_d)\pi(\tau \tilde{\mu}')
d\tilde{\mu}' d\tilde{z}]\nonumber\\
&\leq E_{\Sigma,\tilde{\Sigma}}[\int N(\tilde{z};\Sigma_w^{-1/2}\mu^*,I_d) 
\log \int N(\tilde{z};\tilde{\mu}',I_d)\pi(\tau \tilde{\mu}')
d\tilde{\mu}' d\tilde{z}]
\label{eq:from-thm1}\\
&=E_{\Sigma,\tilde{\Sigma}}[\int N(\tilde{z};\mu^*,\Sigma_w) 
\log \int N(\tilde{z};\tilde{\mu}',\Sigma_w)\pi(\tau \Sigma_w^{-1/2}
\tilde{\mu}')
d\tilde{\mu}' d\tilde{z}]\nonumber
\end{align}
Here, inequality (\ref{eq:from-thm1}) is given by Theorem \ref{thm:reg-1} (ii).

Since every rotation invariant superharmonic function is radially nonincreasing,
$$\pi(\tau \Sigma_w^{-1/2}\tilde{\mu}') \leq \pi(\Sigma^{-1/2}\tilde{\mu}').$$
From this inequality,
\begin{align}
E_{\Sigma,\tilde{\Sigma}}&[\int N(\tilde{z};\mu^*,\Sigma_w) 
\log \int N(\tilde{z};\tilde{\mu}',\Sigma_w)\pi(\tau \Sigma_w^{-1/2}
\tilde{\mu}')
d\tilde{\mu}' d\tilde{z}]\nonumber\\
&\leq E_{\Sigma,\tilde{\Sigma}}[\int N(\tilde{z};\mu^*,\Sigma_w) \log \int 
N(\tilde{z};\tilde{\mu}',\Sigma_w) \pi(\Sigma^{-1/2} \tilde{\mu}')
d\tilde{\mu}' d\tilde{z}] \nonumber\\
&= \varphi_{\pi_{\Sigma}}^w(\mu^*)\nonumber\\
&= \varphi_{\pi_{\Sigma}}^w(\mu)
\end{align}

In particular, if $\pi$ is not constant, inequality
(\ref{eq:from-thm1}) holds strictly.
Therefore, $p_{\Sigma}$ dominates $p_{\rm I}$. \qed

From Lemma \ref{lem:reg-3}, $p_{\Sigma}$ is proved to be minimax.
\begin{corollary}~\\
\label{cor:reg-rot-1}
Assume $d\geq 3$.
Let $p(\Sigma)$ and $\tilde{p}(\tilde{\Sigma})$
be rotation invariant continuous functions.
If $\pi$ is a rotation invariant superharmonic prior,
Bayesian predictive density $p_{\Sigma}(\tilde{y}|y)$
is minimax under $\mathcal{R}_{\rm KL}$.
\end{corollary}

Theorem \ref{thm:reg-rot-2} and Corollary \ref{cor:reg-rot-1} 
can be generalized to the case with a {\it semi}-positive
definite future covariance matrix $\tilde{\Sigma}$.
Let $\tilde{\Sigma}$ be a $d$-dimensional semi-positive matrix whose rank
is $k>0$. Then there is a $d \times k$ matrix $L$ satisfying
$\tilde{\Sigma}=L L^\top$.
Let $\{a_i\}_{i=1}^{d-k}$ be a set of orthogonal normalized
vectors that are orthogonal to each column vector of $L$,
i.e. $L^\top a_i=0$ and $a_i^\top a_j=\delta_{ij}$ for $i,j=1,\dots,d-k$.
Define the Normal distribution with semi-positive definite
covariance matrix by
$$
N_d(y;\mu,\tilde{\Sigma})=
\frac{1}{(2\pi)^{k/2}|L^\top L|^{1/2}}
\exp\left(-\frac{(y-\mu)^\top \tilde{\Sigma}^\dagger (y-\mu)}{2}\right)
\prod_{i=1}^{d-k} \delta(a_i^\top (y-\mu))
$$
where $\tilde{\Sigma}^\dagger$ is Moore-Penrose pseudo-inverse
of $\tilde{\Sigma}$.

From the results of functional analysis, $N_d(y;\mu,\Sigma)$
for any semi-positive definite $\tilde{\Sigma}$ is 
equivalent to $\lim_{\epsilon\rightarrow 0} N_d(y;\mu, \Sigma+\epsilon I_d)$
as a functional on Schwartz functions of $y$.

Using this equivalence and the bounded convergence theorem,
equation (\ref{eq:reg-risk-difference}) is valid
for a semi-definite future covariance matrix
if we define $\Sigma_w := (\Sigma^{-1}+\tilde{\Sigma}^\dagger)^{-1}$.
Because $\tilde{\Sigma}^\dagger \neq 0$, $\tau$
defined by (\ref{eq:tau}) takes value in $(0,1)$.
Therefore, Theorem \ref{thm:reg-rot-2} and Corollary \ref{cor:reg-rot-1}
hold for each semi-definite future covariance matrix $\tilde{\Sigma}$.

\section{Prior distributions depending on the future covariance}
\label{sec:main2}
In this section, we consider prior distributions depending on the future 
covariance matrix. Theorem \ref{thm:new2} below says that every Bayesian 
prediction with an adequately metrized prior dominates that based on the 
uniform prior.
Although the assumption that priors can depend on the future covariance may seem 
strange, this assumption is natural when we consider the linear regression 
problem, as we will see in Section \ref{sec:regression}.

First, we generalize Theorem \ref{thm:reg-1} to the case with 
non-identity covariances.
Let $\mu$ and $z$ be vectors in $\mathbb{R}^d$ and let $\Sigma\in \mathbb{R}^{d\times d}$
be a positive definite matrix. 
	    
Let $\Sigma_1$ and $\Sigma_2$ be positive definite matrices such that 
$\Sigma_1 \preceq \Sigma_2$. An orthogonal matrix $U$ and a diagonal matrix 
$\Lambda$ are given by a diagonalization of $\Sigma_1^{1/2}\Sigma_2^{-1}\Sigma_1^{1/2}$,
i.e. $\Sigma_1^{1/2}\Sigma_2^{-1}\Sigma_1^{1/2}=U^\top\Lambda U$. Let
$A^*:=\Sigma_1^{1/2} U^\top (\Lambda^{-1}-I_d)^{1/2}$.

\begin{proposition}
\label{prop:new1}
If $\pi$ is a prior s.t. $\pi(A^* \mu)$ is a superharmonic function of 
$\mu$,
then
\begin{equation}
\phi_{\pi}(\mu,\Sigma_1)\geq \phi_{\pi}(\mu,\Sigma_2)
\label{eq:new0}
\end{equation}
for any $\mu\in \mathbb{R}^d$.
Inequality (\ref{eq:new0}) becomes strict if $\pi$ is not a constant 
function.
\end{proposition}

The following theorem is a direct result of Proposition \ref{prop:new1}.

\begin{theorem}
\label{thm:new2}
If $\pi(A^* \mu)$ is a superharmonic function of $\mu$, then 
$R_{\rm{KL}}(\pi,\mu)\leq R_{\rm{KL}}(\pi_I,\mu)$.
Furthermore, if $\pi$ is not a constant function, a Bayesian predictive distribution $p_{\pi}$
dominates the one with the uniform prior $\pi_I$.
\end{theorem}

Note that $\pi(A^*\mu)$ can be superharmonic only if ${\rm rank}(\Sigma_2-\Sigma_1)\geq 3$.

\quad\\
\noindent{\bf Proof of Proposition \ref{prop:new1} and Theorem \ref{thm:new2}.}\quad
Assume $0\prec \Sigma_1 \preceq \Sigma_2$ and let 
$\Sigma_1^{1/2}\Sigma_2^{-1}\Sigma_1^{1/2}=U^\top \Lambda U$
be a diagonalization. Then,
\begin{align*}
\phi_\pi(\mu,\Sigma)&=\int \log \left\{\int \pi(\nu)\frac{1}{(2\pi)^{d/2}|\Sigma|^{1/2}}
\exp\left(-\frac{(x-\nu)^\top \Sigma^{-1}(x-\nu)}{2}\right) \diff \nu \right\}\\
&\hspace{2cm}
\frac{1}{(2\pi)^{d/2}|\Sigma|^{1/2}}
\exp\left(-\frac{(x-\mu)^\top \Sigma^{-1}(x-\mu)}{2}\right)\diff x.
\end{align*}

Let $\tilde{x}:=U\Sigma^{-1/2}x$, $\tilde{\mu}=U\Sigma^{-1/2}\mu$, and
$\tilde{\nu}=U\Sigma^{-1/2}\nu$. 
By $|\Sigma_2|^{-1/2} |\Sigma_1|^{1/2}=|\Lambda|^{1/2}$,
\begin{align}
\label{eq:new-1}
\phi_{\pi}(\mu,\Sigma_2)
&=\int \log \left\{\int \pi(\Sigma^{1/2} U^\top \tilde{\nu})
\frac{1}{(2\pi)^{d/2}|\Lambda|^{1/2}}
\exp\left(-\frac{(\tilde{x}-\tilde{\nu})^\top 
\Lambda^{-1}(\tilde{x}-\tilde{\nu})}{2} \right) \diff\tilde{\nu}\right\}\nonumber\\
&\hspace{2cm}
\frac{1}{(2\pi)^{d/2}|\Lambda|^{1/2}}
\exp\left(-\frac{(\tilde{x}-\tilde{\mu})^\top 
\Lambda^{-1}(\tilde{x}-\tilde{\mu})}{2}\right)\diff \tilde{x}\nonumber\\
&=\phi_{\pi(\Sigma_1^{1/2} U^\top \cdot)}(\tilde{\mu},\Lambda^{-1}),
\end{align}
where $\pi(\Sigma_1^{1/2}U^\top\cdot)$ is a prior distribution whose density function 
is represented by $\pi(\Sigma_1^{1/2} U^\top \mu)$ with a prior
density $\pi(\mu)$. 

Putting $\Sigma_2=\Sigma_1$, we get 
\begin{equation}
\phi_{\pi}(\mu,\Sigma_1)=\phi_{\pi(\Sigma_1^{1/2}U^\top\cdot)}(\tilde{\mu},I_d),
\end{equation}
where $I_d$ is the $d$-dimensional identity matrix.

We denote each diagonal component of $\Lambda$ by $\lambda_i$.
Now $0<\lambda_i\leq 1$ for each $i$ since $\Sigma_1 \preceq \Sigma_2$.
Let $a_i(t):= 1+t(\lambda_i^{-1}-1)$ and $A:={\rm diag}(a_i)$.
Then 
\begin{align*}
\phi_\pi&(\mu,\Sigma_2)-\phi_\pi(\mu,\Sigma_1)\\
&=\phi_{\pi(\Sigma_1^{1/2}U^\top\cdot)}(\tilde{\mu},\Lambda^{-1})
-\phi_{\pi(\Sigma_1^{1/2}U^\top\cdot)}(\tilde{\mu},I_d)\\
&=\int_{t=0}^1 \left. \sum_{i=1}^d \frac{\partial a_i(t)}{\partial t}
\frac{\partial}{\partial a_i}
\phi_{\pi(\Sigma_1^{1/2}U^\top\cdot)}(\tilde{\mu},A)
\right|_{a_i(t)} \diff t \\
&=\int_{t=0}^1 \left. \sum_{i=1}^d \frac{\partial \tilde{a}_i(t)}{\partial t}
\frac{\partial}{\partial \tilde{a}_i}
\phi_{\pi(A^* \cdot)}(\hat{\mu},\tilde{A})
\right|_{\tilde{a}_i(t)} \diff t\\
&=\int_{t=0}^1 \left. \sum_{i=1}^d
\frac{\partial}{\partial \tilde{a}_i}
\phi_{\pi(A^* \cdot)}(\hat{\mu},\tilde{A})
\right|_{\tilde{a}_i(t)} \diff t
\end{align*}
where $\tilde{a}_i:=(\lambda_i^{-1}-1)^{-1}a_i$ and 
$\hat{\mu}:=(\Lambda^{-1}-I_d)^{-1/2}\tilde{\mu}$.

By assumption, $\pi(A^* \cdot)$ for 
$A^*=\Sigma_1^{1/2}U^\top (\Lambda^{-1}-I_d)^{1/2}$ is superharmonic.
Now it is sufficient to prove Lemma \ref{lem:new1} iii) below. \qed

\begin{lemma}
\label{lem:new1}
i) $\Sigma^d_{i=1} \frac{\partial}{\partial a_i} N(x;\mu,A)
=\frac{1}{2}\Delta N(x;\mu,A)$.\\
ii) $\int f(x-t) \diff \mu(t)$ is a superharmonic function of $x$ if $f$ 
is a superharmonic function and $\mu$ is a positive measure on $\mathbb{R}^d$.\\
iii) $\sum_{i=1}^d \frac{\partial}{\partial a_i}\phi_\pi(\mu,A)\leq 0$
for any $\mu \in \mathbb{R}^d$, $a_i>0$, and $A={\rm diag}(a_i)$
for each superharmonic prior $\pi$.
\end{lemma}

\quad\\
\noindent{\bf Proof of Lemma \ref{lem:new1}.}\quad
Lemma i) follows from direct calculation.
For a proof of ii), see Problem 1.7.16 of \cite{lehmann_casella1998}.

\begin{align}
\sum_{i=1}^d \frac{\partial}{\partial a_i} \phi_\pi(\mu,A)
&=\sum_{i=1}^d \frac{\partial}{\partial a_i} \int \log \left\{
\int \pi(\nu) N(x;\nu,A)d\nu\}
\right\}N(x;\mu,A)dx \nonumber\\
&=\int \frac{\sum_{i=1}^d \frac{\partial}{\partial a_i}
\int \pi(\nu) N(x;\nu,A)d\nu}{\int \pi(\nu) N(x;\nu,A)d\nu}
N(x;\mu,A)dx \nonumber\\
&+\int  \log \left\{\int \pi(\nu) N(x;\nu,A)d\nu\}
\right\}\sum^d_{i=1} \frac{\partial}{\partial a_i} N(x;\mu,A)dx.
\label{eq:new1}
\end{align}
Now, 
$$
\sum^d_{i=1} \frac{\partial}{\partial a_i} \int \pi(\nu) N(x;\nu,A)\diff \nu
=\frac{1}{2}\Delta \int \pi(\nu)N(x;\nu,A)d\nu \leq 0
$$
from Lemma \ref{lem:new1} i) and ii).
Thus, the first term of the right-hand side of (\ref{eq:new1}) is non-positive.
The second term of the right-hand side of (\ref{eq:new1}) becomes
\begin{align}
\frac{1}{2}\int \log \left\{\int \pi(\nu) N(x;\nu,A)d\nu
\right\}\Delta N(x;\mu,A)dx \nonumber\\
= \frac{1}{2} \int \Delta \log \left\{\int \pi(\nu) N(x;\nu,A)d\nu
\right\} N(x;\mu,A)dx
\label{eq:self-adjoint}
\end{align}
by i) and the self-adjoint property of the Laplacian.
Since the logarithm of a superharmonic function is superharmonic 
(see Problem 1.7.16 of \cite{lehmann_casella1998}),
(\ref{eq:self-adjoint}) is non-positive from ii).
Thus Lemma \ref{lem:new1} iii) is proved. \qed

\begin{example}
\rm
A rescaled Stein prior 
$$\pi_{{\rm S};\Sigma_2-\Sigma_1}(\mu) = \|(\Sigma_2-\Sigma_1)^{-1/2}\mu\|^{-(d-2)}$$
satisfies the condition of Proposition \ref{prop:new1} and Theorem \ref{thm:new2}.
This is because
\begin{align*}
\|(\Sigma_2-\Sigma_1)^{-1/2}\mu\|^{-(d-2)}
&=(\mu^\top \Sigma_1^{-1/2} (\Sigma_1^{-1/2}\Sigma_2\Sigma_1^{-1/2}-I_d)^{-1}
\Sigma_1^{-1/2}\mu)^{-(d-2)/2}\\
&=(\mu^\top \Sigma_1^{-1/2} U^\top (\Lambda^{-1}-I_d)^{-1}U
\Sigma_1^{-1/2}\mu)^{-(d-2)/2}.
\end{align*}
Thus, $\pi_{{\rm S};\Sigma_2-\Sigma_1}(A^* \mu)=\pi_{\rm S}(\mu).$
\end{example}

\section{Application to the Normal linear regression problem}
\label{sec:regression}
In this section, we apply the results
in the previous section to the Normal linear regression problem.

Consider a Normal linear model
\begin{equation}
y = X^{\top} \beta + \epsilon,
\label{eq:reg-normal-linear}
\end{equation}
$$
\epsilon \sim N_p(0, \sigma^2 I_p),
$$
where the target variable $y$ is a $p$ dimensional vector, 
$X$ is a $d \times p$ matrix composed of the explanatory variables,
$\sigma^2>0$ is an unknown variance, and
$\beta$ is an unknown $d$-dimensional vector.
When the rightmost column of $X$ is the constant vector $(1,\dots,1)^{\top}$,
the model (\ref{eq:reg-normal-linear}) is
a model with a intercept, $y = X^{\top} \beta + \beta_0 
+ \epsilon$.

We suppose that a future sample $\tilde{y}$ is generated by
\begin{equation}
\tilde{y} = \tilde{X}^{\top} \beta + \tilde{\epsilon},
\label{eq:reg-normal-linear-future}
\end{equation}
$$
\tilde{\epsilon} \sim N_p(0, \tilde{\sigma}^2 I_p),
$$
where $\tilde{y}$ is a $\tilde{p}$ dimensional vector, $\tilde{X}$ is a 
$d \times \tilde{p}$ matrix,
and $\tilde{\sigma}^2>0$ is an unknown variance.

In the present work, we assume that $p\geq d$ and
$X X^{\top}$ is regular, however neither
$\tilde{p}\geq d$ nor regularity of $\tilde{X}\tilde{X}^\top$
is necessary.

We consider the prediction problem for the linear regression models
(\ref{eq:reg-normal-linear}) and (\ref{eq:reg-normal-linear-future})
with KL risk function
\begin{equation*}
\tilde{R}_{\rm KL}(\beta,\hat{p}_\pi,X,\tilde{X})
:= \int p(y|X;\beta,\sigma^2) D(p(\tilde{y}|\tilde{X};\beta,\tilde{\sigma}^2)
\| p_\pi(\tilde{y}|\tilde{X},y,X;\sigma^2,\tilde{\sigma}^2)) \diff y.
\end{equation*}
and partial Bayesian risk function with prior $p(X)$ and $\tilde{p}(\tilde{X})$:
\begin{equation*}
\tilde{\mathcal{R}}_{\rm KL}(\beta,\hat{p}_\pi)
:=\int p(X) \tilde{p}(\tilde{X}) 
\tilde{R}_{\rm KL}(\beta,\hat{p}_\pi,X,\tilde{X})
\diff X \diff \tilde{X}.
\end{equation*}
Note that we do not assume any prior for $\beta$.

Next, the regression model is reduced to a Normal model
discussed in Section \ref{sec:main_result}.
Let $y_1 :=(X X^{\top})^{-1}X y$ and 
$y_2:= y- X^{\top} (X X^{\top})^{-1}X y$.
Then
\begin{align*}
\frac{1}{(2\pi)^{p/2}}& \exp\Big(-\frac{(X^{\top} \beta - y)^{\top}
(X^{\top} \beta - y)}{2\sigma^2}\Big) \diff y\\
&=
\frac{1}{(2\pi)^{d/2} |\Sigma |^{1/2}}
\exp\Big(-\frac{(y_1- \beta)^{\top} \Sigma^{-1}
(y_1 - \beta)}{2} \Big)
g(y_2;\sigma^2) 
\diff y_1 \diff y_2,
\end{align*}
where 
\begin{equation}
\Sigma := \sigma^2 (X X^{\top})^{-1}
\label{eq:reg-Sigma}
\end{equation}
and $g(y_2;\sigma^2)$ is a density function of $y_2$ that is independent 
of $y_1$ and $\beta$.

When $y$ is given, $y_1$ is a sufficient statistic of $\beta$,
the maximum likelihood estimator, and the least-square estimator of $\beta$.
Thus, the regression model (\ref{eq:reg-normal-linear})
is reduced to a Normal model
\begin{equation}
p(y_1;\beta,\Sigma) = N_d(y_1;\beta, \Sigma).
\label{eq:reg-normal2}
\end{equation}
Similarly, the regression model (\ref{eq:reg-normal-linear-future})
for the future samples is reduced to a Normal model
\begin{equation}
\tilde{p}(\tilde{y}_1;\beta,\tilde{\Sigma}) = N_d(\tilde{y}_1;\beta, \tilde{\Sigma})
\label{eq:reg-normal-future2}
\end{equation}
with semi-positive definite covariance matrix.
Here $\tilde{y}_1 :=(\tilde{X} \tilde{X}^{\top})^\dagger
\tilde{X} \tilde{y}$
and
\begin{equation}
\tilde{\Sigma} := \tilde{\sigma}^2
(\tilde{X} \tilde{X}^{\top})^\dagger.
\label{eq:reg-tilde-Sigma}
\end{equation}

The KL risk of the Bayesian predictive density with a prior 
$\pi(\beta)$ for the regression problem becomes

\begin{align}
\tilde{R}_{\rm KL}(p_\pi,\beta)
&= \int p(y|X;\beta,\sigma^2) 
D(p(\tilde{y}|\tilde{X};\beta,\tilde{\sigma}^2)\|
p_\pi(\tilde{y}|\tilde{X},y,X)) \diff y\nonumber\\
&= \int p(y|X;\beta,\sigma^2)
\int N_d(\tilde{y}_1;\beta,\tilde{\Sigma}) 
g(\tilde{y}_2;\tilde{\sigma}^2)\nonumber\\
&
\log  \frac{N_d(\tilde{y}_1;\beta,\tilde{\Sigma})
g(\tilde{y}_2;\tilde{\sigma}^2)}
{\displaystyle \frac{\int N_d(\tilde{y}_1;\beta,\tilde{\Sigma}) 
g(\tilde{y}_2;\tilde{\sigma}^2)
N_d(y_1;\beta,\Sigma) g(y_2;\sigma^2) \pi (\beta) \diff \beta}
{\int N_d(y_1;\beta,\Sigma) g(y_2;\sigma^2)\pi(\beta) \diff \beta}}
\diff \tilde{y}_1 \diff \tilde{y}_2
\diff y
\nonumber\\
&= \int N_d(y_1;\beta,\Sigma) D(N_d(\tilde{y}_1;\beta,\tilde{\Sigma})
\| q_\pi(\tilde{y}_1|y_1)) \diff y_1\nonumber\\
&= R_{\rm KL}(q_\pi,\beta),
\end{align}
where
$$
q_\pi(\tilde{y}_1|y_1):= \frac{\int N_d(\tilde{y}_1;\beta,\tilde{\Sigma})
N_d(y_1;\beta,\Sigma) \pi(\beta) \diff \beta}
{\int N_d(y_1;\beta,\Sigma) \pi(\beta) \diff \beta}.
$$

As a result, the prediction problem for the regression model
(\ref{eq:reg-normal-linear}) and (\ref{eq:reg-normal-linear-future})
is reduced to a prediction problem (\ref{eq:reg-normal2}) and
(\ref{eq:reg-normal-future2}).
Using the result in Section \ref{sec:main_result}, 
we construct a Bayesian prediction for the Normal regression problem.

Define $\Sigma$, $\tilde{\Sigma}$, and $\Sigma_w$ by (\ref{eq:reg-Sigma}), 
(\ref{eq:reg-tilde-Sigma}), and $\Sigma_w = (\Sigma^{-1}+\tilde{\Sigma}^\dagger)^{-1}$,
respectively, then the following theorem and corollary hold.

\begin{theorem}\quad
\label{thm:reg-reg-main}

Let $\pi_{\Sigma}(\beta)=\pi(\Sigma^{-1/2} \beta)$.
Let $p(X)$ and $\tilde{p}(\tilde{X})$
be rotation invariant continuous functions.

(i) If $\pi$ is a non-constant rotation invariant superharmonic function,
then the Bayesian predictive density $p_{\Sigma}$ with a prior
$\pi_{\Sigma}$ dominates
$p_{\rm I}$ with the uniform prior $\pi_{\rm I}$ under the risk $\tilde{\mathcal{R}}_{\rm KL}$.

(ii) If $\pi$ is a rotation invariant superharmonic function,
then $p_{\Sigma}$ is minimax under the KL risk $\tilde{\mathcal{R}}_{\rm KL}$.
\end{theorem}

\Proof 
If $p(X)$ and $\tilde{p}(\tilde{X})$ are rotation invariant,
then the distributions of $\Sigma = \sigma^2 (X X^{\top})^{-1}$ and
$\Sigma_w = (\sigma^{-2} (X X^{\top})+\tilde{\sigma}^{-2} (\tilde{X} \tilde{X}^{\top}))^{-1}$
are also rotation invariant.

From Theorem \ref{thm:reg-rot-2} and Corollary \ref{cor:reg-rot-1},
the theorem is derived directly. \qed

The assumption of rotation invariance of $p(x)$ and $p(\tilde{x})$
is sometimes not realistic. If we consider priors depending on the future 
explanatory variables, we can construct a Bayesian prediction dominating
the one with the uniform prior and, therefore, being a minimax prediction.

Define an orthogonal matrix $U$ and a diagonal matrix $\Lambda$ by a 
diagonalization of $\Sigma_w^{1/2} \Sigma^{-1} \Sigma_w^{1/2}$, i.e.
$\Sigma_w^{1/2} \Sigma^{-1} \Sigma_w^{1/2}=U^\top \Lambda U$.
Let $A^*:= \Sigma_w^{1/2} U^\top (\Lambda^{-1}-I_d)^{1/2}$.
Then the following theorem is a direct consequence of Theorem \ref{thm:new2}.

\begin{theorem}
\label{thm:reg-new1}
(i) If $\pi(A^* \beta)$ is superharmonic w.r.t. $\beta$ and $\pi$ is non-constant,
then the Bayesian prediction based on the prior $\pi$ dominates that based on the uniform prior.\\
(ii)If $\pi(A^* \beta)$ is superharmonic, then the Bayesian prediction 
based on the prior $\pi$ is minimax.
\end{theorem}

Note that $\pi(A^*\beta)$ can be superharmonic only if the number of 
the future samples is more than two.

\section{Experimental results}
We show several experimental results on
the Bayesian prediction with shrinkage priors
for regression problems.

\begin{figure}[t]
\begin{center}
    \includegraphics[width=8cm]{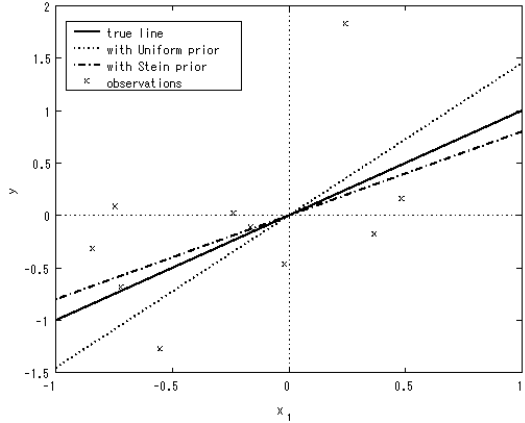}
  \caption{An example of the Bayesian prediction based on the uniform prior 
and a rescaled Stein prior for the Normal regression model without an intercept term.
}
\label{fig:reg-example-nobias}
\end{center}
\end{figure}

\begin{figure}[tbh]
\begin{center}
    \includegraphics[width=8cm]{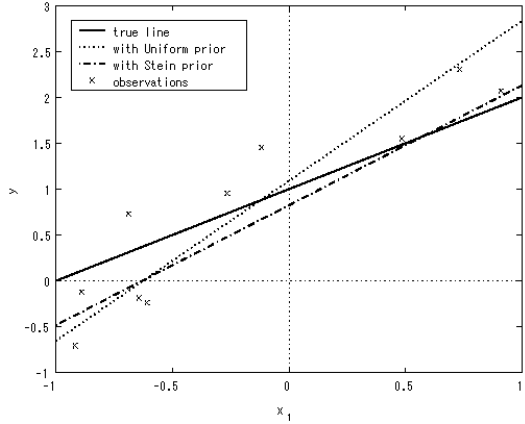}
\caption{An example of the Bayesian prediction based on the uniform prior 
and a rescaled Stein prior for the Normal regression model with an intercept term
$\beta_0=1$.
}
\label{fig:reg-example-bias}
\end{center}
\end{figure}

Figures \ref{fig:reg-example-nobias} and \ref{fig:reg-example-bias}
are examples of the regression problem.
We consider the five dimensional Normal regression models, without
an intercept term (Figure \ref{fig:reg-example-nobias}) and with an intercept term
(Figure \ref{fig:reg-example-bias}).
We set the true parameter $\beta=(1,0,\dots,0)\in \mathbb{R}^5$.
An explanatory variable $X$ is sampled from the uniform distribution
$U([-1,1]^{5\times 10})$ and
corresponding target variable $y$ is sampled from $N_{10}(X^\top \beta,I_{10})$.
The target variable $\tilde{y}$ for each explanatory variable
$\tilde{x}=(\tilde{x}_1,0,\dots,0)$ where
$\tilde{x}_1\in [0,2]$ is predicted by the Bayesian predictive density
based on the uniform prior $\pi_{\rm I}$ and that based on a rescaled Stein prior
$\pi_{S;\Sigma}$ where $\Sigma = XX^\top$.

Two lines in Figures \ref{fig:reg-example-nobias} and \ref{fig:reg-example-bias}
are $y=\hat{\beta}_\pi^\top \tilde{x}$ for $\pi_{\rm I}$ and $\pi_{S;\Sigma}$, respectively,
where $\hat{\beta}_\pi$ is the posterior mean
with prior $\pi$.
In both figures,
the slope of the line with rescaled Stein prior is smaller than the one with
the uniform prior because the slope parameter $\beta$ is shrunk to $\beta=0$.
Moreover in Figure \ref{fig:reg-example-bias},
the intercept parameter is also shrunk.

\begin{figure}[tbh]
\begin{center}
    \includegraphics[width=8cm]{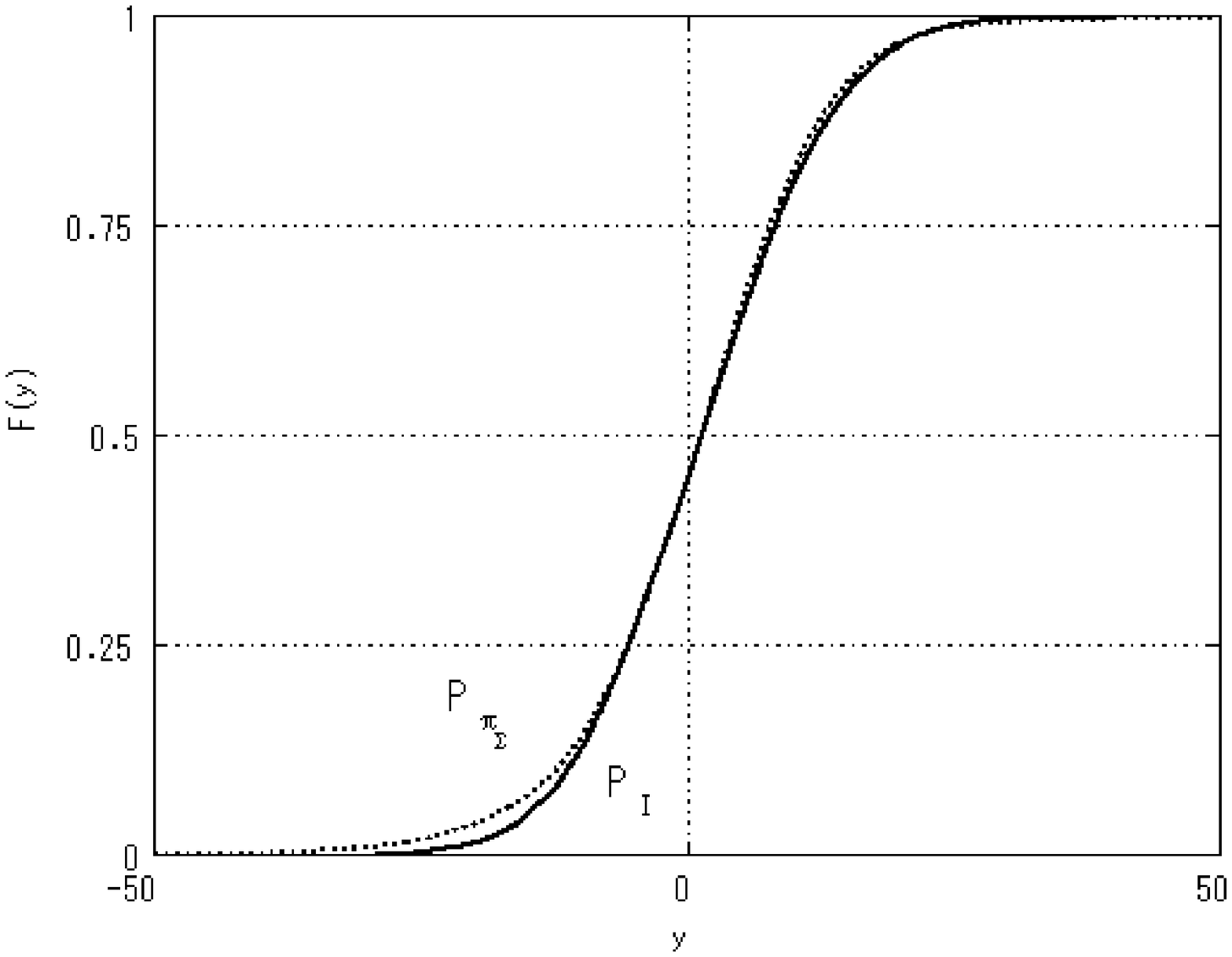}
  \caption{Distribution functions of $p_{\rm I}(\tilde{y}|\tilde{x},y,X)$ and
$p_{{\rm S};\Sigma}(\tilde{y}|\tilde{x},y,X)$ where 
$\beta=\tilde{x}=e_1:=(1,0,\dots,0)\in \mathbb{R}^d$,
$X$ is a sample from $U([-1,1]^{d\times p})$, $y$ is a sample from
$N_p(y;X^\top \beta, 10I_p)$, and $\tilde{p}(\tilde{y}|\tilde{x})={\rm N}(\tilde{y};
\tilde{x}^\top \beta,10)$. We generate $10^4$ samples of $\tilde{y}$ from 
each predictive distribution. 
Sample means of $P_{I}$ and $P_{{\rm S};\Sigma}$ are $1.3134$ and $0.6898$, respectively.
}
\label{fig:reg-histgram}
\end{center}
\end{figure}

Figure \ref{fig:reg-histgram} shows
the distribution functions of the predictive density 
$p_{\rm I}(\tilde{y}|\tilde{x},y,X)$ with $\pi_{\rm I}$ and
$p_{{\rm S};\Sigma}(\tilde{y}|\tilde{x},y,X)$ with $\pi_{{\rm S};\Sigma}$,
respectively, for $\beta=\tilde{x}=e_1:=(1,0,\dots,0)\in \mathbb{R}^d$.

\begin{figure}[tbh]
\begin{center}
    \includegraphics[width=8cm]{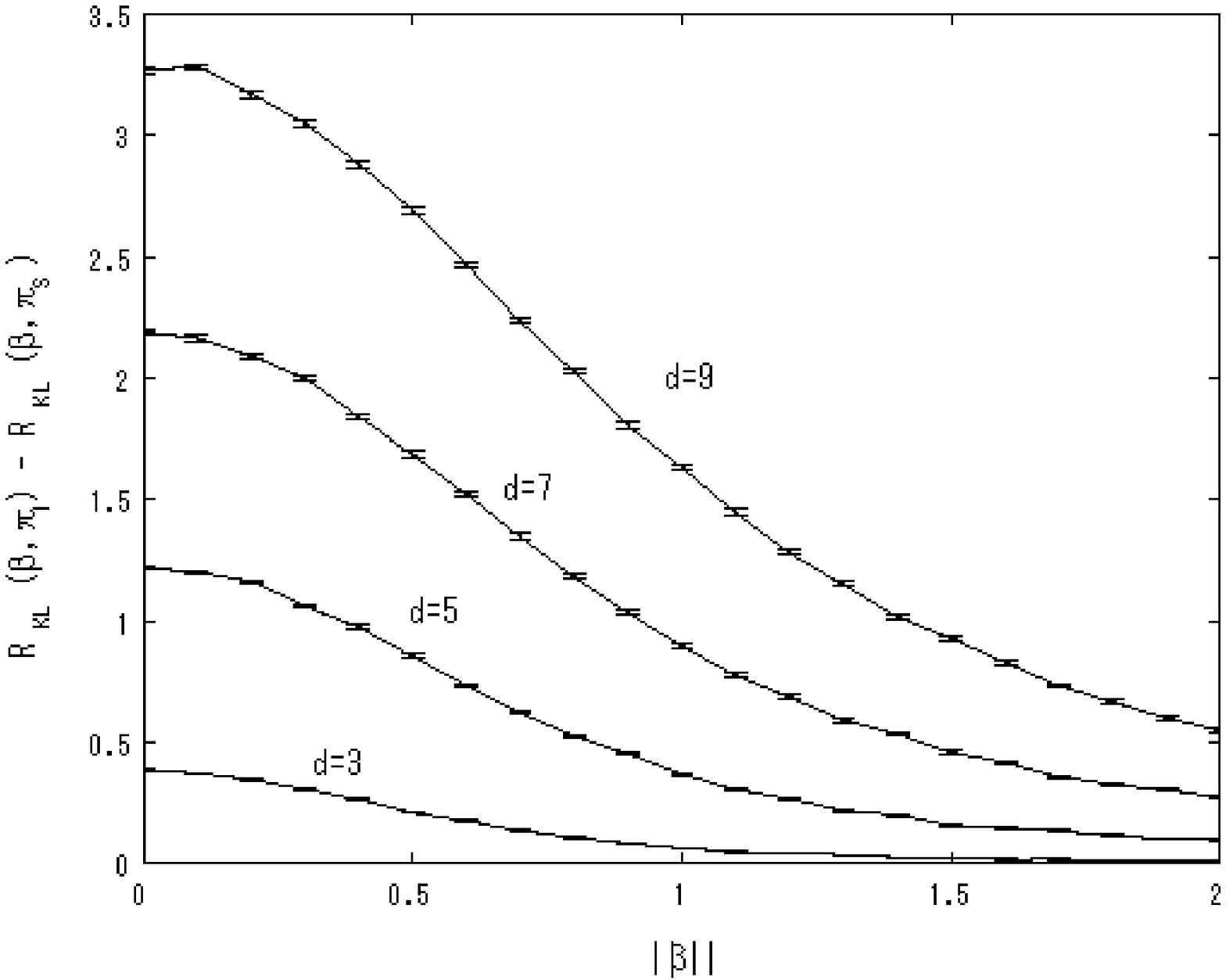}
  \caption{The risk difference of $p_I$ and $p_{S;\Sigma}$
for $d=3,5,7,9$ and $\|\beta\|\ in [0,2].$ 
We generate $10^4$ independent samples of $X$ and $\tilde{X}$ from $N_{10}(0,I_{10})$.
Each line in the figure represents the sample mean of risk difference 
$R_{\rm KL}(\beta,p_{\rm I})-R_{\rm KL}(\beta,p_{{\rm S};\Sigma})$.
Each error bar represents the standard deviation.}
\label{fig:reg-dim-nobias}
\end{center}
\end{figure}

\begin{figure}[tbh]
\begin{center}
    \includegraphics[width=8cm]{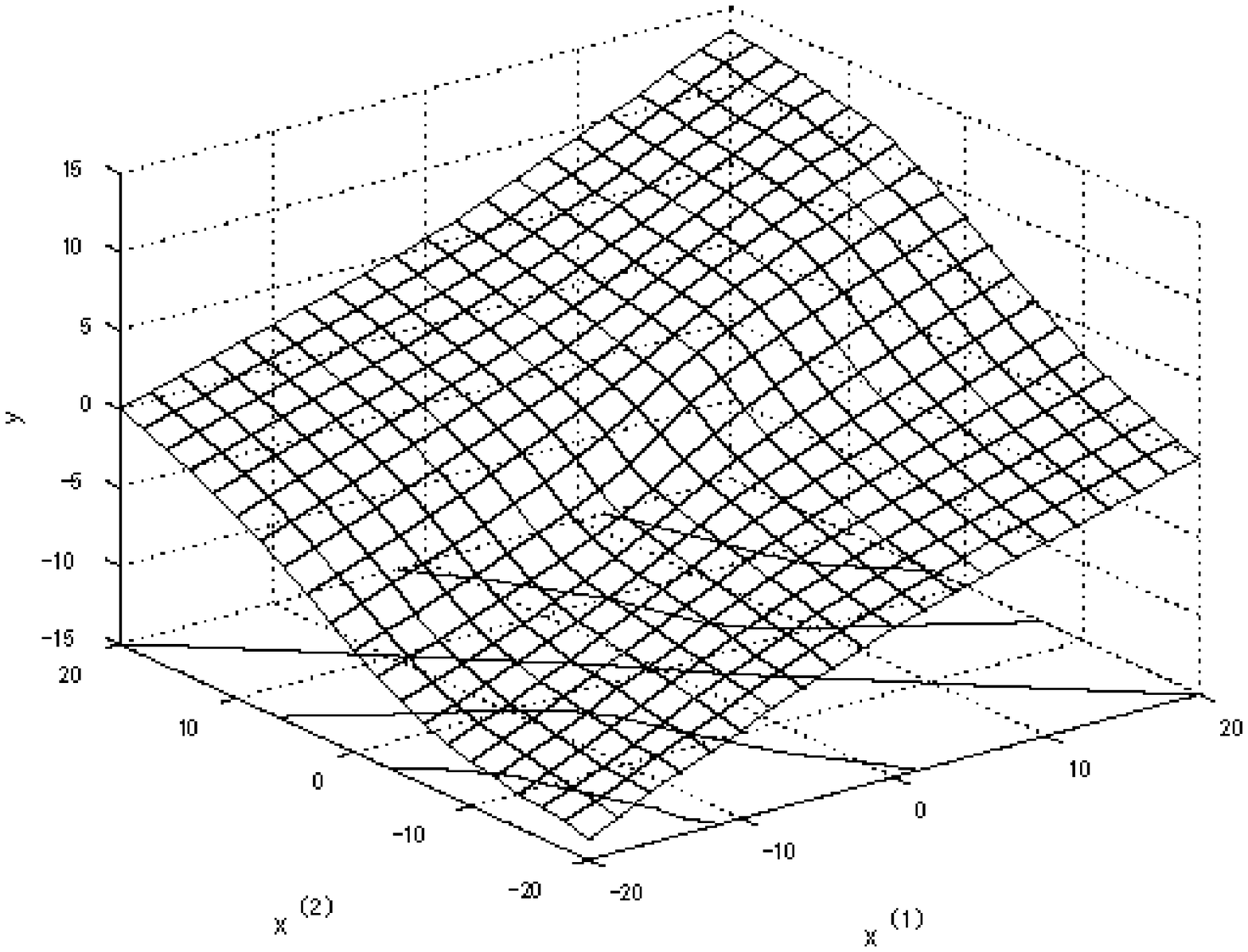}
\caption{An example of Bayesian prediction whose prior depends on the 
explanatory variables of future samples.}
\label{fig:reg-new1}
\end{center}
\end{figure}

Next, we show an example of Bayesian prediction whose prior depends on 
the explanatory variables of future samples. We set 
$x_1=(\sqrt{3}/2,1/2,0)^\top$, $x_2=(\sqrt{3}/2,-1/2,0)^\top$, $x_3=(0,0,1)^\top$,
$y_1=\sqrt{3}/2+1/2$, $y_2=\sqrt{3}/2-1/2$ and $y_3=0$.
Figure \ref{fig:reg-new1} is a graph of 
$E_{\pi_{S;A^*}}[\tilde{y}|\tilde{x},y,x]$ for each value of 
$\tilde{x}=(\tilde{x}^{(1)},\tilde{x}^{(2)},0)\in \mathbb{R}\times 
\mathbb{R} \times \{0\}$ with the rescaled Stein prior $\pi_{\Sigma;A^*}$.
Here, the Bayesian estimation based on the uniform prior corresponds to 
the MLE $\hat{\beta}=(1,1,0)$, i.e. $y=x^{(1)}+x^{(2)}$.

We can see that the amount of shrinkage by the Bayesian prediction increases as
the direction of $\tilde{x}$ becomes closer to $x^{(1)}$
than $x^{(2)}$, i.e. $\tilde{x}^\top e_1$ becomes larger than $\tilde{x}^\top 
e_2$. This fact is intuitively explained as follows:
when explanatory variables of training samples are closer to $x^{(1)}$,
$\tilde{x}$ whose direction is close to $x^{(1)}$ has more information 
than ones whose direction is close to $x^{(2)}$. Thus $\tilde{x}$ close 
to $x^{(1)}$ need not be shrunk.

Figure \ref{fig:reg-dim-nobias} shows the risk functions of $p_I$ and $p_\Sigma$
for $d=3,5,7,9$ and $\|\beta\|\in [0,2]$.
The model has no intercept term.
Here we assume that the columns of 
$X$ and $\tilde{X}$ are independently sampled from
$N_{10}(0,I_{10})$.

\begin{figure}[tbh]
\begin{center}
    \includegraphics[width=8cm]{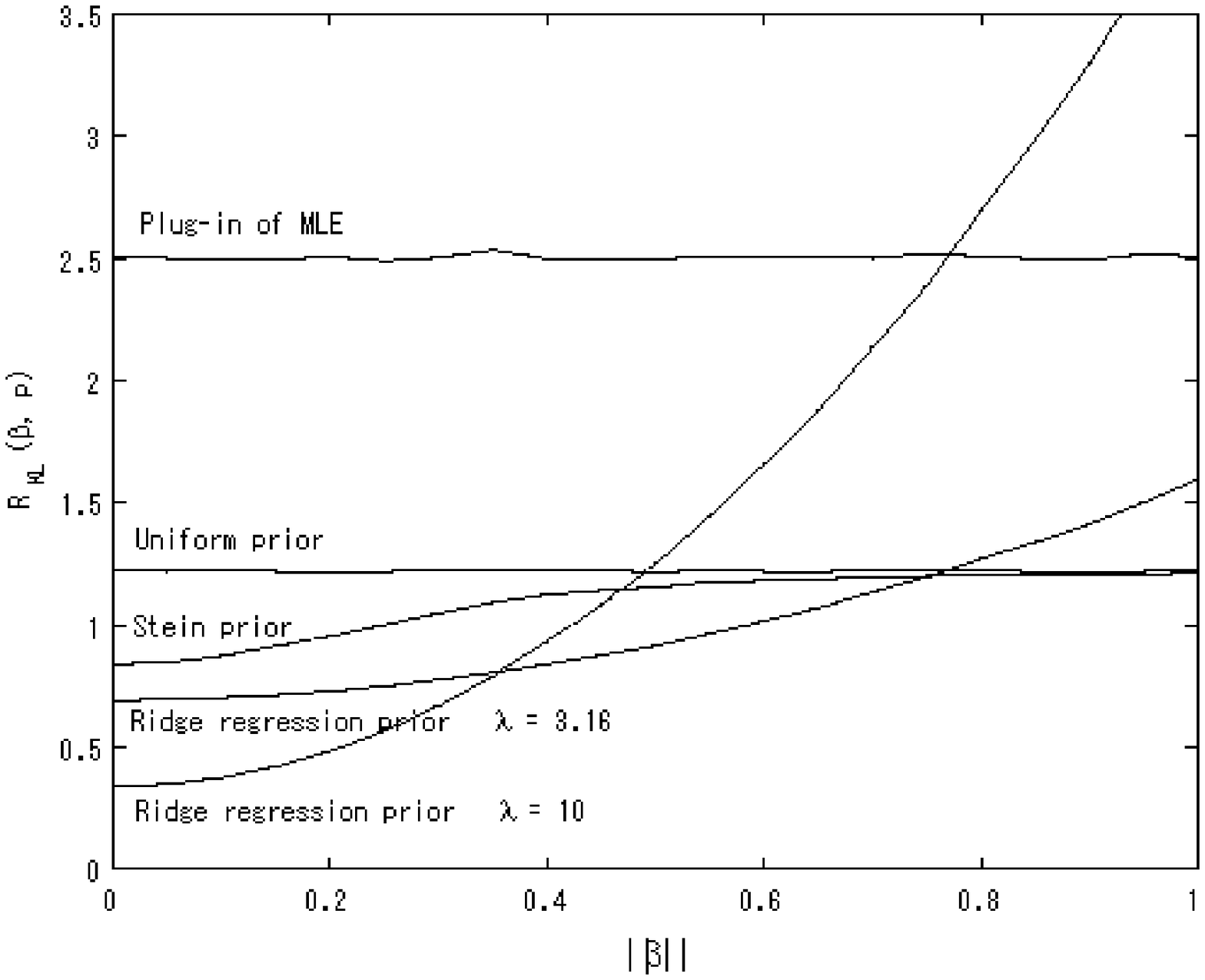}
  \caption{Comparison of the risk values by five predictive densities:
the Bayesian predictive density based on $p_I$ and $p_{\Sigma}$,
the ridge regression prior with regularization parameters $\lambda = 10$ 
and $\lambda = \sqrt{10}=3.16$,
and the plug-in density of the MLE. The model is five dimensional 
and has no intercept term.
We generate $10^4$ independent samples of $X$ and $\tilde{X}$ from $N_{10}(0,I_{10})$.
Each line in the figure represents the sample mean of the risk $R_{\rm KL}(\beta,\hat{p})$
for the predictive density $\hat{p}$.}
\label{fig:reg-risks}
\end{center}
\end{figure}

Figure \ref{fig:reg-risks} compares
five predictive densities: the Bayesian predictive density based on $p_I$ and $p_{\pi_\Sigma}$,
the ridge regression prior with regularization parameters $\lambda \in 
\{\sqrt{10}, 10\}$,
and the plug-in density of MLE.

The ridge regression prior is
$$
\pi_{RR}(\beta;\lambda)=\frac{\lambda^{d/2}}{(2\pi)^{d/2}}
\exp\Big(-\lambda\frac{\|\beta\|^2}{2}\Big)
$$
with a regularization parameter $\lambda >0$.
We note that the posterior mean with the ridge regression prior
is equivalent to the ridge regression estimator
$$
\hat{\beta}_{RR}= (XX^\top+ \lambda I)^{-1}Xy.
$$

When $\|\beta\|$ is close to $0$, the center of shrinkage,
the risk based on the ridge regression prior $\pi_{RR}$
becomes smaller than that based on $\pi_\Sigma$.
However, when $\|\beta\|$ increases,
the prediction with $\pi_{RR}$ becomes worse than
the one with $\pi_{\rm I}$ and even worse than 
the plug-in distribution of the MLE. 

\section{Conclusions and discussions}

In this paper, we considered the multivariate Normal model with an unknown mean and
a known covariance.
The covariance matrix can be changed after the first sampling.
We assumed rotation invariant priors of
the covariance matrix and the future covariance matrix.
We showed that the shrinkage predictive density with 
the rescaled rotation invariant superharmonic priors 
is minimax under the Kullback-Leibler risk.
Moreover, if the prior is not constant,
Bayesian predicitive density based on the prior
dominates the one with the uniform prior.

In this case, the rescaled priors are independent of
the covariance matrix of future samples.
Therefore, we can calculate the posterior distribution
and the mean of the predictive distribution
(i.e. the posterior mean and the Bayesian estimate for quadratic loss) 
based on some of the rescaled Stein priors without knowledge of
future covariance.
Since the predictive density with the uniform prior is minimax,
the one with each rescaled Stein prior is also minimax.

Next we considered Bayesian predictions whose prior can depend on the 
future covariance. In this case, we proved that the Bayesian prediction 
based on a rescaled superharmonic prior dominates the one with the uniform 
prior without assuming the rotation invariance.

Applying these results to the prediction
of response variables in the Normal regression model, 
we show that there exists the prior distribution such that
the corresponding Bayesian predictive density
dominates that based on the uniform prior.
Since the prior distribution is independent of future explanatory variables,
both the posterior distribution and the mean of the predictive distribution
are independent of the future explanatory variables.

The robustness of some shrinkage methods as Stein estimators
has been studied (see, for example, the bibliography in \cite{robert2001}).
The Stein effect has robustness in the sense that
it depends on the loss function rather than the true distribution of
the observations.
Our result shows that the Stein effect has
robustness with respect to the covariance of the true distribution of
the future observations.

As the dimension of the model becomes large, 
the risk improvement by the shrinkage with the rescaled Stein prior
$\pi_\Sigma$ increases as in Figure
\ref{fig:reg-dim-nobias}.
An important example of the high dimensional model
is the kernel methods (see \cite{ESL}).
As noted in \cite{intro_to_svm},
the feature space of kernel methods 
is a kernel reproducing Hilbert space 
whose dimension is as large as the sample size.
Therefore Bayesian prediction based on shrinkage priors 
could be efficient for kernel methods. This is a future problem.

\section{Acknowledgment}
The authors appreciate Mr. Vu, Vincent Q. for precious comments on an 
earlier version of this paper.

\appendix
\section{Finiteness of the marginal distribution}
Here, we prove finiteness of the marginal distribution $m_\pi(\mu,\Sigma)$.
\begin{lemma}
\label{lem:finite-marginal}
If $\pi$ is a superharmonic prior density function, the marginal distribution $m_\pi(x,\Sigma)$ is 
finite for every vector $x\in \mathbb{R}^d$ and positive 
definite matrix $\Sigma\in \mathbb{R}^{d\times d}$.
\end{lemma}

\Proof
Fix a vector $x\in \mathbb{R}^{d}$
From the definition of superharmonic functions, $\pi\not\equiv \infty$.
Thus, $\exists x_0\in \mathbb{R}^d$ s.t. $\pi(x_0)<\infty$.
If we set $\tilde{\pi}(\mu):=\pi(\mu + x_0)$, then $\tilde{\pi}$ is 
superharmonic and $\tilde{\pi}(0)<\infty$.

Let $\lambda_{\rm max}$ be the maximal eigenvalue of $\Sigma$ and 
$r_0:=\| x+x_0 \|$, then
\begin{align}
m_\pi(x,\Sigma)
&\leq \int \exp \left(-\frac{\|x+x_0-\mu \|^2}{2 
\lambda_{\rm max}}\right)\tilde{\pi}(\mu)\diff{\mu}\nonumber\\
&\leq \int_{\|\mu\|\leq 2r_0} 
\exp \left(-\frac{\|x+x_0-\mu \|^2}{2 
\lambda_{\rm max}}\right)\tilde{\pi}(\mu)\diff{\mu}
+ \int_{\|\mu\|> 
2r_0}\exp\left(-\frac{\|\mu\|^2}{8\lambda_{\rm max}}\right)
\tilde{\pi}(\mu)\diff{\mu}
\label{eq:ap-1}
\end{align}
The first term of the right-hand side of (\ref{eq:ap-1}) is finite
because the integral of a superharmonic function over a compact subspace
of $\mathbb{R}^d$ is finite (see Theorem 4.10 of \cite{helms-1969}).

The second term is also finite because
$$
\sum_{n=2}^\infty \int_{n r_0<\|\mu\|\leq (n+1)r_0}
\exp\left(-\frac{\|\mu\|^2}{8\lambda_{\rm max}}\right)
\tilde{\pi}(\mu) \diff \mu
\leq C\sum_{n=2}^\infty \exp\left(-\frac{(nr_0)^2}{8\lambda_{\rm max}}\right)
\tilde{\pi}(0) \{(n+1)r_0\}^d
$$
for a positive constant $C$. Here we used a fact 
$\int_{\|\mu\|<r}\tilde{\pi}(\mu) \diff \mu < C \tilde{\pi}(0)r^d$
by Theorem 4.9 of \cite{helms-1969}.
Therefore, $m_\pi(x,\Sigma)<\infty$. \qed

From this lemma, we see the assumption $m_\pi(z,v I_d)<\infty$ in 
Theorem \ref{thm:reg-1} (ii) is redundant.


\begin{thebibliography}{13}
\expandafter\ifx\csname natexlab\endcsname\relax\def\natexlab#1{#1}\fi

\bibitem[{Corcuera \& Giummol{\'e}(2000)}]{corcuera_giummole2000}
\textsc{Corcuera, J.~M.} \& \textsc{Giummol{\'e}, F.} (2000).
\newblock First-order optimal prediction densities.
\newblock In \textit{Applications of differential geometry to econometrics},
  P.~Marriott \& M.~Salmon, eds. Cambridge: Cambridge University Press, pp.
  214--229.

\bibitem[{Cristianini \& Shawe-Taylor(2000)}]{intro_to_svm}
\textsc{Cristianini, N.} \& \textsc{Shawe-Taylor, J.} (2000).
\newblock \textit{An Introduction to Support Vector Machines}.
\newblock Cambridge University Press.

\bibitem[{George et~al.(2006)George, Liang \& Xu}]{george_etc2006}
\textsc{George, E.~I.}, \textsc{Liang, F.} \& \textsc{Xu, X.} (2006).
\newblock Improved minimax prediction under {K}ullback-{L}eibler loss.
\newblock \textit{Annals of Statistics} \textbf{34}, 78--91.

\bibitem[{Hastie et~al.(2001)Hastie, Tibshirani \& Friedman}]{ESL}
\textsc{Hastie, T.}, \textsc{Tibshirani, R.} \& \textsc{Friedman, J.} (2001).
\newblock \textit{The elements of statistical learning -- Data mining,
  inference, and prediction}.
\newblock springer series in statistics. New York: Springer.

\bibitem[{Helms(1969)}]{helms-1969}
\textsc{Helms, L.~L.} (1969).
\newblock \textit{Introduction to Potential Theory}.
\newblock New York: Wiley-Interscience.

\bibitem[{Komaki(1996)}]{komaki1996}
\textsc{Komaki, F.} (1996).
\newblock On asymptotic properties of predictive distributions.
\newblock \textit{Biometrika} \textbf{83}, 299--313.

\bibitem[{Komaki(2001)}]{komaki2001}
\textsc{Komaki, F.} (2001).
\newblock A shrinkage predictive distribution for multivariate {N}ormal
  observables.
\newblock \textit{Biometrika} \textbf{88}, 859--864.

\bibitem[{Komaki(2006)}]{komaki-2006AS}
\textsc{Komaki, F.} (2006).
\newblock Shrinkage priors for bayesian prediction.
\newblock \textit{Annals of Statistics} \textbf{34}, 808--819.

\bibitem[{Lehmann \& Casella(1998)}]{lehmann_casella1998}
\textsc{Lehmann, E.~L.} \& \textsc{Casella, G.} (1998).
\newblock \textit{Theory of point estimation}.
\newblock New York: Springer, 2nd ed.

\bibitem[{Liang \& Barron(2004)}]{liang-barron-2004IEEEIT}
\textsc{Liang, F.} \& \textsc{Barron, A.} (2004).
\newblock Exact minimax strategies for predictive density estimation, data
  compression, and model selection.
\newblock \textit{Ieee Transactions on Information Theory} \textbf{50},
  2708--2726.

\bibitem[{Murray(1977)}]{murray1977}
\textsc{Murray, G.~D.} (1977).
\newblock A note on the estimation of probability density functions.
\newblock \textit{Biometrika} \textbf{64}, 150--152.

\bibitem[{Ng(1980)}]{ng1980}
\textsc{Ng, V.~M.} (1980).
\newblock On the estimation of parametric density functions.
\newblock \textit{Biometrika} \textbf{67}, 505--506.

\bibitem[{Robert(2001)}]{robert2001}
\textsc{Robert, C.~P.} (2001).
\newblock \textit{The {B}ayesian Choice}.
\newblock New York: Springer-Verlag, 2nd ed.

\end{thebibliography}

\end{document}